\documentclass{amsart}
\usepackage{amsmath}
\usepackage{amssymb}
\usepackage{amsthm}
\usepackage{color}
\usepackage{parskip}
\usepackage{import}
\usepackage[noadjust]{cite}
\usepackage{mathtools}
\usepackage[final]{fixme}
\usepackage[final]{changes}

\makeatletter
\def\thm@space@setup{%
  \thm@preskip=\parskip \thm@postskip=0pt
}
\makeatother

\theoremstyle{plain}
\newtheorem{theorem}{Theorem}[section]
\newtheorem{lemma}[theorem]{Lemma}
\newtheorem{proposition}[theorem]{Proposition}
\newtheorem{corollary}[theorem]{Corollary}
\newtheorem{conjecture}[theorem]{Conjecture}

\theoremstyle{definition}
\newtheorem{remark}{Remark}[section]

\newtheorem{example}{Example}[section]
\newtheorem{definition}{Definition}[section]
\DeclareMathOperator{\per}{per}
\DeclareMathOperator{\ind}{ind}
\DeclareMathOperator{\Br}{Br}
\DeclareMathOperator{\im}{im}
\DeclareMathOperator{\br}{br}
\DeclareMathOperator{\pr}{pr}
\DeclareMathOperator{\Fp}{\mathbb{F}_p}
\DeclareMathOperator{\Hom}{Hom}

\DeclareMathOperator{\chara}{char}
\DeclareMathOperator{\F}{\mathbb{F}}

\DeclareMathOperator{\Z}{\mathbb{Z}}
\DeclareMathOperator{\id}{id}
\DeclareMathOperator{\tx}{\tilde{x}}
\DeclareMathOperator{\sab}{\mathcal{S}_{a,b}}

\newcommand{\brac}[1]{\left( #1\right)}

\newcommand{\blue}[1]{\textcolor{blue}{#1}}
\newcommand{\red}[1]{\textcolor{red}{#1}}

\usepackage{hyperref}
\hypersetup{
    colorlinks=true,
    urlcolor=black,
    citecolor=blue,   
    linkcolor=blue,
}

\title{Brauer $p$-dimension of complete discretely valued fields}
\author{Nivedita Bhaskhar}
\address{Department of Mathematics
University of California at Los Angeles \\
Los Angeles, CA 90095-1555}
\email{nbhaskh@math.ucla.edu}
\author{Bastian Haase}
\address{Department of Mathematics \& Computer Science, Emory University, 400 Dowman Drive
NE, Atlanta, GA 30322, USA}
\email{bhaase@emory.edu}

\begin{document}
\maketitle      
                       

\begin{abstract}
  Let $K$ be a complete discretely valued field of characteristic 0 with residue field $\kappa$ of characteristic $p$. Let $n=[\kappa:\kappa^p]$ be the $p$-rank
of $\kappa$. It was proved in \cite{Parimala2013} that the Brauer $p$-dimension of $K$ lies between $n/2$ and $2n$. For $n\leq 3$,
 we improve the upper bound to $n+1$ and provide examples to show that our bound is sharp. For $n \leq 2$, we also improve the lower bound
to $n$.
For general $n$, we construct a family of fields $K_n$ with residue fields of $p$-rank $n$, such that $K_n$ admits a central simple algebra $D_n$ of index $p^{n+1}$.
Our sharp lower bounds for $n\leq 2$ and upper bounds for $n\leq 3$ in combination with the nature of these examples motivate us to conjecture that the Brauer $p$-dimension of such fields always lies between $n$ and $n+1$. \par
\deleted[id=Bastian]{We exploit Kato's filtration of the $p$-torsion of the Brauer group as in \cite{Parimala2013} along with a deeper analysis of the symbol length of $k_2(\kappa)$.}
\end{abstract}

\section{Introduction}

\let\thefootnote\relax\footnote{\textbf{2010 Mathematics Subject Classification:} 16K50, 11R58. \par
\textbf{Key words and phrases:} Brauer group, Brauer $p$-dimension, complete discretely valued field, $p$-rank, Milnor $k$-groups, Kato's filtration.}

To every central simple algebra $A$ over a field $F$, one can associate two numerical Brauer class invariants. One is the \textit{index}, denoted $\ind(A)$ which is the degree of the unique (up to isomorphism) central division algebra Brauer equivalent to $A$. The other is the \textit{period},  denoted $\per(A)$, which is the order of the Brauer class of $A$ in the Brauer group $\Br(F)$. It is well known from the classical theory of central simple algebras that $\per(A)|\ind(A)$ and in fact that $\per(A)$ and $\ind(A)$ also have the same set of prime factors.

Thus one can define the \textit{Brauer dimension} of a field $F$, denoted $\Br\dim(F)$, to be the least number $n$ such that $\ind(A)|\per(A)^n$ for every central simple algebra $A/E$ for
any finite extension $E/F$. If no such $n$ exists, then $\Br\dim(F)$ is set to $\infty$.

The Brauer dimension of separably closed fields and finite fields is clearly $0$ since these fields have trivial Brauer group. From local class field theory, we can see that the Brauer dimension of a local field is $1$ . The same result holds for a global field as a consequence of the famous theorem of Albert-Hasse-Brauer-Noether. The \textit{period-index questions} which involve bounding the Brauer dimensions of arbitrary fields form an integral part of this research area.

Let $p$ be a prime. One can define the \textit{Brauer $p$-dimension} of a field $F$, denoted $\Br_p\dim(F)$, to be the least number $n$ such that $\ind(A)|\per(A)^n$ for every central simple algebra $A/E$ of period a power of $p$ for any finite extension $E/F$. As before, if no such $n$ exists, then $\Br_p\dim(F)$ is set to $\infty$. 

Let $F$ be the function field of a $p$-adic curve. One of the long standing period-index questions was  whether the Brauer $\ell$-dimension of such an $F$ is at most $2$ (or even finite!).  For $\ell\neq p$, this was affirmatively answered by Saltman in the late 90s (\cite{Saltman1997}). In the so-called \textit{bad characteristic case}, i.e. $\ell=p$, this was again answered affirmatively by the more recent work of Parimala and Suresh (\cite{Parimala2013}) using patching techniques of Harbater-Hartmann-Krashen.

Recall that the $p$-rank of a characteristic $p$ field $\kappa$ is $n$ if $[\kappa : \kappa^p] = p^n$. Thus a characteristic $p$ field $\kappa$ is perfect if and only if its $p$-rank is $0$. In fact, in the same paper, Parimala and Suresh also investigate the Brauer $p$-dimension of function fields of curves over complete discretely valued fields whose residue fields are not necessarily perfect  and obtain the following more general result: 

\begin{theorem}[\cite{Parimala2013},Thm 3]
Let $K$ be a complete discretely valued field  with residue field $\kappa$. Suppose that $\chara(K) = 0$ and $\chara(\kappa) = p > 0$. Let $F$ be the function field of a curve over $K$. If the $p$-rank of $\kappa$ is $n$, then $\Br_p\dim(F)\leq 2n+2$.
\end{theorem}

A crucial step in the proof involves estimating the Brauer $p$-dimensions of complete discretely valued fields with characteristic $p$ residue fields of $p$-rank $n+1$. This step is essential for applying the patching techniques of HHK.

Using the tools provided by Kato's isomorphisms between the filtrations of the Milnor $K$-group modulo $p$ and the $p$-torsion of the Brauer group of $F$, they have the following estimate:

\begin{theorem}[\cite{Parimala2013},Thm 2.7]
Let $K$ be a complete discretely valued field of characteristic 0 with residue field $\kappa$. Suppose that $\chara(\kappa) = p > 0$ and the $p$-rank of $\kappa$ is $n$, then $\lfloor{n/2}\rfloor \leq \Br_p\dim(K)\leq 2n$.
\end{theorem}

While the proof is the consequence of some very subtle manipulations of Kato's filtrations, the bounds in the above theorem do not appear to be optimal. In this paper, we investigate the Brauer $p$-dimensions of complete discretely valued fields with low $p$-rank residue fields and find  better upper bounds. More precisely, we obtain the following theorem:

\begin{theorem}[Theorem \ref{BN-mainthm1}]
Let $K$ be a complete discretely valued field of characteristic 0 with residue field $\kappa$. Suppose that $\chara(\kappa) = p > 0$ and the $p$-rank of $\kappa$ is $n$ where $n=0,1, 2$ or $3$. Then $\Br_p\dim(K)\leq n+1$.  \par
For $n < 3$, we have $n \leq \Br_p \dim(K)$ and for $n=3$ we have $2 \leq \Br_p\dim(K)$. \end{theorem}

 Our approach, while also utilizing Kato's filtrations, differs from that adopted by Parimala and Suresh in that we rely on bounding the symbol length of the second Milnor $K$ group modulo $p$, denoted by $k_2\brac{\kappa}$, in a concrete manner. This analysis in turn relies on the machinery of differentials in characteristic $p$ as developed by Cartier and the isomorphism between the group of logarithmic differentials $\nu(2)_{\kappa}$ and $k_2(\kappa)$.

We also show that these bounds are optimal when $n \leq 2$ by providing relevant examples:

\begin{theorem}[Theorem \ref{exampleleq2}]
Let $p$ be a prime and let $n, i$ be integers such that $0\leq n\leq 2$ and $n\leq i\leq n+1$. 
Then there exists a characteristic $0$ field $K$ with residue field $k$ of characteristic $p$ and $p$-rank $n$ whose Brauer-$p$-dimension is $i$.
\end{theorem}

This paper is structured as follows.

The next section is devoted to recalling some facts from the theory of differentials and Cartier's theorem which gives useful criteria to identify boundaries and cycles. In the same section, we also recall a very useful but slightly convoluted filtration of the module of differentials of a characteristic $p$-field $\kappa$. In the third section, the proof of the surjectivity of the map $k_2\brac{\kappa}\rightarrow \nu(2)_{\kappa}$ is re-examined carefully for low $p$-ranks to understand the shape of any element  in \replaced[id=Bastian]{$\nu(2)_k$}{this group} in a more concrete fashion (c.f. Theorem \ref{rep}). In the fourth section, we use this theorem to understand the Brauer $p$-dimension of fields as in Theorem \ref{BN-mainthm1}, which emboldens us to make the following conjecture: 

\begin{conjecture}[Conjecture \ref{conjBH}]
Let $K$ be a complete discretely valued field with residue field $\kappa$. Suppose that $\chara(\kappa) = p > 0$ and that the $p$-rank of $\kappa$ is $n$. Then $n \leq \Br_p\dim(K)\leq n+1$. 
\end{conjecture}

\replaced[id=Bastian]{A positive solution to the above conjecture  would be immensely useful in making the Brauer $p$-dimension
 bounds more precise in the corresponding function field cases.}{We hope that in the future this conjecture might be proven to hold.
This in turn will be useful in making the Brauer $p$-dimension bounds more precise in the corresponding function field cases. }

In the fifth and final section, we set forth examples realizing possible Brauer $p$-dimensions in the low $p$-rank cases investigated. Additionally, for each $n\geq 1$, we construct complete discretely valued fields with residue fields of $p$-rank $n$ such that they admit central simple algebras of index $p^{n+1}$. This shows that the optimal upper bound for the Brauer $p$-dimension of such fields cannot be less than $n+1$.

\section*{Acknowledgements}
We would like to thank Professors R. Parimala and V. Suresh for many valuable discussions and critical comments. 
The authors are  partially supported by National Science Foundation grants DMS-1401319 and DMS-1463882.
\section{Preliminaries}

\subsection{Brauer group of complete discretely valued fields}
Let $R$ be a complete discretely valued ring with field of fractions $K$ and residue field $k$. For any finite field extension $L$ of $K$, we can uniquely extend the valuation of $K$ to $L$. Let the valuation ring of $L$ be $S$ and the residue field be $l$. We denote the \emph{ramification index} of $L$ over $K$ by $e$. It is the unique integer $e$ satisfying $\pi_K=u\pi_L^e$ where $\pi_K,\pi_L$ are parameters of $K$ and $L$ respectively and $u$ is a unit in $S$. Let $f$ denote $[l:k]$. If $n=[L:K]$, then it is well-known that $n=ef$.

If $e=1$ and $l$ is a separable extension of $k$, then we say that the field extension $L/K$ is \emph{unramified}. Otherwise, we say that it is \emph{ramified}.

\added[id=Bastian]{Recall the notion of ramification of central simple algebras over $K$.}
\begin{definition}
 Let $A$ be a central simple algebra over $K$ and let $[A]$ denote its class in the Brauer group of $K$. Then, we say that $A$ is \emph{unramified} if
$[A]$ is in the image of the natural injection $\Br(R)\hookrightarrow \Br(K)$. 
Otherwise, we say that $A$ is \emph{ramified}.
\end{definition}

\begin{lemma}
Let $K$ be a complete discretely valued field and let $L$ be a nontrivial, finite, unramified field extension and let $\pi_K$ be a 
parameter of $K$. Then, $\pi_K$ is not a norm from $L$.
\label{norm}
\end{lemma}
\begin{proof}
\replaced[id=Bastian]{This follows from the fact that $\nu_K(N(x))=f\nu_L(x)$ for $x \in L$ where $f=[l:k]$.}{
Since $L/K$ is unramified, $\pi_K$ is a parameter for $L$ also.  Let now $x=u \pi_K^n$ be an arbitrary element in $L$ with $u \in S^*$ and $n \in \mathbb{Z}$. As the norm of a unit in $S$ is a unit in $R$, we see that $N_{L/K}(x) = N_{L/K}(u)\pi_K^{n[L:K]}$. Since $n[L:K]\neq 1$, clearly $\pi_K$ cannot be a norm from $L$.}
\end{proof}

\begin{lemma}
Let $R'$ be a cyclic \'{e}tale algebra over $R$ with generator $\sigma$ and let $u \in R^{*}$. Then, 
$\left( R' ,\sigma,u\right)$ is an Azumaya algebra.
\label{unramified}
\end{lemma}
\begin{proof}
  This follows readily from Proposition 1.2b) in chapter IV of \cite{Milne1980}.
\end{proof}

Let $F$ be a field of characteristic $p$. For $a,b\in F^*$, we can define the \emph{symbol $p$-algebra} 
\begin{align*}
  [a,b):=F[x,y \;|\; x^p-x=a, y^p=b, yx=xy+y]
\end{align*}
which is always central simple over $F$ of period dividing $p$. 

\begin{proposition}[\cite{Gille2006},Corollary 4.7.5]
  Let $a,b$ in  $F$. Then, $[a,b)$ is split  iff $b$ is a norm from $F[t]/(t^p-t-a)$.
\label{normsplit}
\end{proposition}

In particular, $[a,b)$ splits if $b$ is a $p$-th power or
$x^p-x-a$ is split over $F$.

 Recall that there are natural isomorphisms  (cf. Corollary 2.5 in \cite{Vishne2000})
\begin{align*}
  [a_1,b)\otimes[a_2,b)&\simeq [a_1+a_2,b)\otimes M_p(F) \\
  [a,b_1)\otimes [a,b_2)& \simeq [a,b_1b_2)\otimes M_p(F).
\end{align*}

Let $K$ now be a compelete discretely valued field of characteristic $0$ with residue field $k$ of characteristic $p$.
Let $a,b \in K^*$ such that $x^p-x-a$ is irreducible. Then note that $L=K[x]/(x^p-x-a)$ is a Galois extension over $K$ as we can lift all roots from $k[x]/(x^p-x-a)$ by Nakayama's Lemma.
Let $\sigma$ denote a lift of the automorphism of $k[x]/(x^p-x-\overline{a})$ defined via $x \mapsto x+1$. Then, $\sigma$ is a generator of the Galois group of $L$ over $K$. We then define 
\begin{align*}
  [a,b):=K[x,y \; | \; x^p-x=a,y^p=b,yx=\sigma(x)y].
\end{align*}
Let $\overline{a},\overline{b}$ denote the residues of $a,b$. Then, we have that $[\overline{a},\overline{b})=[a,b) \in \Br(K)$ where
$[\overline{a},\overline{b})$ is identified with its image under $\Br(k)\xrightarrow{\sim} \Br(R) \hookrightarrow \Br(K)$.

We will use these algebras when constructing examples in section \ref{examples}.

\subsection{Differentials in characteristic $p$}
In this section, we will collect some well-known results on differentials of fields of characteristic $p$. The main reference for most of the results here
are \cite{Colliot-Thelene1997} and \cite{Gille2006}. \par

Let $k$ be a field of characteristic $p$ and $F$ a subfield containing $k^p$.
Let $\Omega^1_{k/F}$ denote the module of differentials of $k$ relative to $F$ and let $\Omega^q_{k/F}=\bigwedge_{i=1}^q \Omega^1_{k/F}$.
These modules form a complex with the usual differential map $d^q: \Omega^q_{k/F}\rightarrow \Omega^{q+1}_{k/F}$. \deleted[id=Bastian]{We will often omit $q$ when it is clear from
context.} Let $Z^q_{k/F}$ denote the kernel of $d^q$ and let $B^q_{k/F}$ denote the image of $d^{q-1}$. 
\subsubsection{Cartier's Theorem}
 Define group homomorphisms
\begin{align*}
  \gamma_q\colon \Omega^q_{k/F}& \longrightarrow Z_{k/F}^q/B^q_{k/F},\\
 x \frac{d y_1}{y_1} \wedge\ldots \wedge \frac{d y_q}{y_q}& \mapsto x^p \frac{d y_1}{y_1} \wedge
\ldots \wedge \frac{d y_q}{y_q} 
\end{align*}
for non-negative integers $q$.\deleted[id=Bastian]{ We will often omit $q$ if it is clear from context.}
Cartier proved that $\gamma_q$ is an isomorphism. Let us denote the inverse  by $C_q$.
Note that it satisfies the property $C_q(x^p \frac{d y_1}{y_1} \wedge
\ldots \wedge \frac{d y_q}{y_q} )=x \frac{d y_1}{y_1} \wedge
\ldots \wedge \frac{d y_q}{y_q} $.\par
Recall that we say that a differential $w \in \Omega^1_{k/F}$ is \emph{exact} if $w=da$ and \emph{logarithmic} if $w=\frac{da}{a}$ for some $a \in k^*$.
The maps $\gamma_1$ and $C_1$ allow one to give an equivalent description of these properties.
\begin{theorem}[Cartier]
A differential $w \in \Omega^1_{k/F}$ is exact if and only if $dw=0$ and $C(w)=0$.
  A differential $w \in Z^1_{k/F}$ is logarithmic if and only if 
$\gamma(w)=w \in Z^1_{k/F} / B^1_{k/F}$. \par
\label{cartier}
\end{theorem}
\begin{proof}
   Compare Corollary 1.2.3 in \cite{Colliot-Thelene1997}. Note that our definition of $C_q$ slightly differs from the definition given there.
\end{proof}
This result leads us to the definition of the group $\nu(q)_{k/F}$, which generalizes the notion of being logarithmic. 
\begin{definition}
  The group $\nu(q)_{k/F}$ is the kernel of the map 
  \begin{align*}
     \gamma_q-\id_{\Omega^q_{k/F}}: \Omega^q_{k/F} \rightarrow \Omega_{k/F}^q/B^q_{k/F}.
  \end{align*}
  In the case $F=k^p$, we will just write $\nu(q)_{k}$ for the corresponding group of differentials.
\end{definition} 
\begin{remark}
  Note that the morphism $\gamma-1$ is a generalization of the Artin-Schreier map.
  Furthermore, $\nu(q)_{k/F}$ is functorial in $F$.
\end{remark}
\subsubsection{A filtration on $\Omega_k$} \label{filtration}
We will now introduce a filtration on $\Omega_k$ depending on the choice of a $p$-basis. Suppose that $[k:F]=p^n$ and fix an ordered $p$-basis $\{b_1, \ldots , b_n\}$ of $k$ over $F$. 
Enroute to defining the filtration, we will collect several useful lemmata.
\begin{lemma}
  Let $k_0=F(b_{1},\ldots,b_{l})$ for some $1\leq l \leq n$ and let $q$ be a positive integer.
 Then, the kernel of the natural projection $\Omega^{q}_{k}\rightarrow \Omega_{k/k_0}^q$ is, as a  $k$ vector space, generated by elements of the form 
 \begin{align*}
   \frac{d b_{s(1)}}{b_{s(1)}} \wedge\ldots \wedge \frac{d b_{s(q)}}{b_{s(q)}}
 \end{align*}
where $s:\{1,\ldots,q\}\rightarrow \{1,\ldots,n\}$ is strictly increasing with $s(1)\leq l$.
\label{kernel}
\end{lemma} 
\begin{proof}
  Immediate from the definition of $\Omega^q_{k/k_0}$.
\end{proof}
We can now start decomposing our differentials.
 For any mapping $\mu: \{1,\ldots,n\} \rightarrow \{0,1, \ldots, p-1 \}$, we use the notation
$b^{\mu}=b_{1}^{\mu(1)}\cdots b_{n}^{\mu(n)}$. This gives us 
\begin{align*}
  k&=\bigoplus_{\mu} F \cdot b^{\mu} 
\end{align*}
For any positive integer $q$, let $S_q$ be the set of strictly increasing mappings $s:\{1,\ldots,q\}\rightarrow \{1,\ldots,n\}$.
Then, for $s \in S_q$, we denote 
\begin{align*}
  w_s = \frac{d b_{s(1)}}{b_{s(1)}} \wedge\ldots \wedge \frac{d b_{s(q)}}{b_{s(q)}} \in \Omega^q_{k/F}
\end{align*}
and therefore obtain 
\begin{align*}
  \Omega^q_{k/F}=\bigoplus_{s\in S_q}k\cdot w_s.
\end{align*}
For fixed $\mu$, we denote 
\begin{align*}
  \Omega^{q}_{k/F}(\mu)= \bigoplus_{s \in S_q} F \cdot b^{\mu} w_s \subset \Omega^{q}_{k/F}
\end{align*}
so that we obtain 
\begin{align}
  \Omega^q_{k/F}=\bigoplus_{\mu}\Omega^q_{k/F}(\mu).
  \label{omega_decom}
\end{align}
From now on, when we speak of the \emph{zero component} of an element $w \in \Omega^q_{k/F}$, we mean the element $w_0$
in $\Omega^q_{k/F}(0)$ in the decomposition of $w$ with respect to (\ref{omega_decom}).
Finally, let us  equip $S_q$ with the lexicographical order. For $s \in S_q$, define
\begin{align*}
  \Omega^q_{k/F, < s} =\bigoplus_{t \in S_q, t<s} k \cdot w_{t}
\end{align*}
for fixed $s \in S_q$. This gives us our  filtration of $\Omega^q_{k/F}$.
\begin{lemma}[\cite{Colliot-Thelene1997}, Lemma 3.1]
  Let $k$ and $F$ be as above. Then, 
  \begin{itemize}
    \item We have $d \left( \Omega^q_{k/F}(0) \right) =0$  and
          $d \left(\Omega^q_{k/F}(\mu) \right)\subset \Omega^{q+1}_{k/F}(\mu) $  for all $q \geq0$
    \item For $\mu \neq 0$, the following sequence of $F$ vector spaces is exact 
      \begin{align*}
        0 \xrightarrow{} Fb^{\mu} \xrightarrow{d} \Omega^{1}_{k/F}(\mu) \xrightarrow{d} \Omega^2_{k/F}(\mu)
        \xrightarrow{d} \ldots \xrightarrow{d} \Omega^n_{k/F}(\mu) \xrightarrow{} 0.
      \end{align*}

  \end{itemize}
\label{exact}
\end{lemma} 
\begin{remark}
Note that the lemma above implies that for any element $w \in \ker(d)$, we have that $w \in \im(d)$ if and only if the zero component of $w$ is trivial.
We will use this fact repeatedly in the next section.
\end{remark}
Let us now end this section with some lemmata that will come in handy shortly.
\begin{lemma}
Let $\{a,b\} \subset \kappa$ be a $p$-basis of a field $\kappa$ of characteristic $p$. Then, $c\in \kappa$ is a $p$-th  power if and only if neither $\{a,c\}$ nor $\{b,c\}$ are a $p$-basis. 
  \label{ppower}
\end{lemma} 
\begin{proof}
  Let $c=\sum_{0\leq i,j <p} \lambda_{ij}a^ib^j$ for $\lambda_{ij}\in \kappa^p$. Note that $c$ is a $p$-th power if and only if $\lambda_{ij}=0$ whenever
$i>0$ or $j>0$. Also, recall that $da$ and $db$ form a basis of $\Omega^1_k$. Now, $\{a,c\}$ being $p$-dependent is equivalent to 
$da\wedge dc=0$. This in turn is equivalent to $\lambda_{ij}=0$ for $j>0$. Similarly, $\{b,c\}$ being $p$-dependent is equivalent to
$\lambda_{ij}=0$ for $i>0$. Hence $c\in \kappa^p$ iff $\lambda_{ij}=0$ for $i>0$ or $j>0$ which happens iff both $\{a,c\}$ and $\{b,c\}$ are $p$-dependent. 
\end{proof}
\begin{lemma}[\cite{Colliot-Thelene1997}, Lemma 3.2]
  Let $k$ be a field of characteristic $p$. Suppose that $k$ does not admit any extensions of degree prime to $p$.
Let $E=k(b)$ for $b^p \in k$ but $b \notin k$. Then, for any $k-linear$ mapping $g:E\rightarrow k$, there is 
$c \in E^{*}$ such that $g(c^i)=0$ for $i=1,\ldots,p-1$. \par
\added[id=Bastian]{If $g(1) \neq0$, then $c \notin k$.}
  \label{super-zero}
\end{lemma}

\begin{lemma}
 Let $x \in k \setminus{F}$ and $w\in \Omega_{k/F}$. If $w \wedge \frac{dx}{x}=0 \in \Omega_{k/F}$,
then  $w=0 \in \Omega_{k/F(x)}$.
\label{cycletest}
\end{lemma} 
\begin{proof}
We will only  prove the case when $w \in \Omega^2_{k/F}$ to simplify notation and since this is the only case
we will use.
  Fix a $p$-basis $\{x=x_1,x_2,\ldots, x_n\}$ of $k$ over $k_0$. Then, when we express
$w$ in terms of this basis, the assumption implies that we obtain 
\begin{align*}
  w= \sum_{i=2}^n \lambda_i \frac{d x_1}{x_1} \wedge \frac{d x_i}{x_i}
\end{align*} 
for some $\lambda_i \in k$. Note that $w$ cannot have a summand of the form $ \lambda \frac{d x_j}{x_j} \wedge \frac{d x_i}{x_i}$ with $1<i<j$  as this would contradict
the assumption $w\wedge \frac{dx}{x}=0$. In view of this equality, the results follows from Lemma \ref{kernel}.
\end{proof}

\begin{lemma}
  Let $\{b_1,\ldots,b_s\}$ be a $p$-basis of $k$ over $k_0$ and let $\{y_1,\ldots,y_r\}$ be $p$-independent over $k_0$. Let $a\in k$ be such that 
  \begin{align}
    b_1^ia \frac{d y_1}{y_1} \wedge\ldots \wedge \frac{d y_r}{y_r} \in d\Omega^{r-1}_{k/k_0}
\label{eq:decomp0}
  \end{align}
holds for $1\leq i \leq p-1$. Then, $a=e_0+e_1$ where $e_0 \in k_0$ and $e_1 \in V$ where
$V$ is the $k_0(b_1)$ vector-space generated by $b_2^{j_2}\cdots b_s^{j_s}$ for $0\leq j_k \leq p-1$
and $\sum_k j_k >0$.
\label{decomp}
\end{lemma} 
\begin{proof}
 Since $\{b_i\}$ forms a  $p$-basis, we can write 
  \begin{align}
    a= e_0 +\alpha_1b_1 + \ldots +\alpha_{p-1}b_1^{p-1}+e_1
    \label{eq:decomp}
  \end{align}
for $e_0,\alpha_i \in k_0$ and $e_1 \in V$. We wish to show $\alpha_i=0$.
For $1\leq i \leq p-1$, multiplying equation (\ref{eq:decomp}) by $b_1^i$ yields
\begin{align*}
   b_1^i a= e_0b_1^i +\alpha_1b_1^{i+1} + \ldots + \alpha_{p-1}b_1^{p+i-1}+e_1b_1^{i}
  \end{align*}
By assumption and Lemma \ref{exact}, the zero component  of expression (\ref{eq:decomp0})
is zero. This precisely  means $\alpha_{j}=0$ for $p-1\geq j \geq 1$.
\end{proof} 
\begin{lemma}
  Let $\{b_1,\ldots, b_s\}$ be a $p$-basis of $k$ over $k_0$ and let $0\neq a=e_0+e_1$ with
$e_0 \in k_0$ and $e_1\in V$ where $V$ denotes the $k_0(b_1)$ vector space generated
by $ b_2^{j_2}\cdots b_s^{j_s}$ for $0\leq j_k \leq p-1$ and $\sum_k j_k>0$. Assume furthermore  that
\begin{align}
  (a^p-a)\frac{d b_1}{b_1} \wedge \frac{d b_2}{b_2} \in d \Omega^1_{k/k_0}.
\label{eq:symbol}
\end{align}
Then,  there is an element $y \in k \setminus k_0(b_1)$ such that in $\Omega^2_{k/k_0}$,
\begin{align*}
  a \frac{d b_1}{b_1} \wedge \frac{d b_2}{b_2}=\frac{d b_1}{b_1} \wedge \frac{d y}{y}.
\end{align*}

\label{symbol}
\end{lemma} 
\begin{proof}
We want to apply Cartier's theorem, i.e. Theorem \ref{cartier} to $a\frac{db_2}{b_2}$.
  Since $d \circ d=0$, the assumption implies that $d(a^p-a)\wedge \frac{d b_1}{b_1} \wedge \frac{d b_2}{b_2} =0 \in \Omega^3_{k/k_0}$. As $d(a^p)=0$,  we can derive $d(a)\wedge \frac{d b_1}{b_1} \wedge \frac{d b_2}{b_2}=0$. By Lemma \ref{cycletest}, this implies $da\wedge \frac{db_2}{b_2}=0$
in $\Omega^2_{k/k_{0}(b_1)}$. So, $a\frac{db_2}{b_2}$ and therefore also $(a^p-a)\frac{db_2}{b_2}$
are cycles in $\Omega^1_{k/k_0(b_1)}$. By means of Lemma \ref{exact} and the assumption (\ref{eq:symbol}),
we can derive that $a^p=e_0$. This implies that $a^p-a=-e_1 \in V$ holds. Note that by definition of $V$, $e_1$ has no zero component in
$k$ over $k_0(b_1)$. Thus, using Lemma \ref{exact} again, we can conclude that $(a^p-a)\frac{db_2}{b_2}$ is a boundary in $\Omega^1_{k/k_0(b_1)}$. It follows that the assumption of Theorem \ref{cartier} is satisfied. Hence, there is $y \in k$ such that $a\frac{db_2}{b_2}=\frac{dy}{y} $ in $\Omega^1_{k/k_0(b_1)}$. By means of Lemma \ref{kernel}, we can conclude 
\begin{align*}
  a\frac{db_2}{b_2}=\frac{dy}{y}+e\frac{db_1}{b_1}
\end{align*} 
in $\Omega_{k/k_0}^1$ for some $e \in k$. Therefore, we get the desired equality 
\begin{align*}
  a \frac{d b_1}{b_1} \wedge \frac{d b_2}{b_2} = \frac{d b_1}{b_1} \wedge \frac{d y}{y}
\end{align*}
in $\Omega^2_{k/k_0}$. It is clear that $y \in k \setminus k_0(b_1)$, for otherwise we would
have $\frac{d b_1}{b_1} \wedge \frac{d y}{y}=0$ in $\Omega^2_{k/k_0}$ which would contradict the assumption that $a\neq0$.
\end{proof}

\begin{lemma} \label{prank} 
Let $l/k$ be a finite or separable field extension of fields of characteristic $p$ and let $n$ be the $p$-rank of $k$.
Then, the $p$-rank of $l$ is also $n$. Furthermore, if the extension is separable, then a $p$-basis of $k$
is also a $p$-basis of $l$.
\label{prank}
\end{lemma}
\begin{proof}
A proof of the first part for finite extensions can be found in \cite{Bourbaki1974a}, A.V.135, Corollary 3.
To prove the second statement, let us assume that $l/k$ is finite and separable. Thus, we have
$l=k(\alpha)$ for some $\alpha \in l$.
Let $u_1,\ldots, u_n$ be a $p$-basis of $k$ over $k^p$. The result now follows from
\[l^p(u_1,\ldots,u_n) = k^p(\alpha^p)(u_1,\ldots,u_n)= k(\alpha^p) = l.\] \par
Let now $l/k$ be an infinite separable extension and let $u_1, \ldots, u_n$ denote a $p$-basis of $k$. 
If $\{u_i\}$  were $p$-dependent over $l$, then they would
be $p$-dependent over some finite extension of $k$. As this contradicts the finite case, we conclude that they are $p$-independent. 
Let now $x \in l$. Then, $k(x)/k$ is finite separable and $u_1,\ldots, u_n$ forms a $p$-basis of $k(x)$. Hence, 
$x \in k(x)^p(u_1, \ldots, u_n) \subset l^p(u_1, \ldots, u_n)$ and the result follows.

\end{proof}

\section{Symbol Length in $k_2(k)$}
Throughout this section, let $k$ be a field of characteristic $p$ and write $F=k^p$. Furthermore, assume that $k$ does not admit any extensions of
degree prime to $p$. Let $n$ denote the $p$-rank of $k$, i.e. $[k:k^p]=p^n$.
Recall the isomorphism (\cite{Kato1982},\cite{Bloch1986})
\begin{align*}
  k_2(k) &\rightarrow \nu(2)_k\\
  \{x,y \}& \mapsto \frac{d x}{x} \wedge \frac{d y}{y}
\end{align*}
 where $k_q(k)=K_q(k)/p$ denotes the $q$-th Milnor $K$ group modulo $p$.

We will now analyze the proof of surjectivity for $p$-rank 2 and 3 carefully following \cite{Colliot-Thelene1997}/\cite{Gille2006} to get more
information about elements in $\nu(2)_k$. We will start with three lemmata.

In the following, let $\{ b_1, \ldots, b_n\}$ denote a $p$-basis of $k$.
\begin{lemma}\label{one}
Let $w=a \frac{d b_1}{b_1} \wedge \frac{d b_2}{b_2}$ be an element in $\nu(2)_k$ with $a \neq 0$. Then, there are $z_1 \in F(b_1) \setminus F$ and
$z_2 \in F(b_1,b_2) \setminus F(b_1)$ such that 
\begin{align*}
  w= \frac{d z_1}{z_1} \wedge \frac{d z_2}{z_2}.
\end{align*}
\end{lemma}
\begin{proof}
Let us define $k_2=F(b_1,b_2)$ and $k_1=F(b_1)$. 
We  first claim that $a \in k_2$. This is because of the following: Since $a \frac{d b_1}{b_1} \wedge \frac{d b_2}{b_2} \in \nu(2)_k$, we have $da\wedge \frac{d b_1}{b_1} \wedge \frac{d b_2}{b_2}=0$. If $a \notin k_2$, the set $\{b_1,b_2,a\}$ would be $p$-independent  over $F$ and therefore $da \wedge \frac{d b_1}{b_1} \wedge \frac{d b_2}{b_2}$ would certainly not be zero. 
Hence, $a \in k_2$, so that our assumptions carry over to $\Omega_{k_2/F}$.\par
We want to apply Lemma \ref{super-zero} on the mapping 
\begin{align*}
  g \colon k_1 a \frac{d b_1}{b_1} \wedge \frac{d b_2}{b_2} \subset \Omega^2_{k_2/F}
\xrightarrow{\pr} \Omega^2_{k_2/F} /d \Omega^1_{k_2/F}.
\end{align*}
Note that the dimension of the $F$ vector space $\Omega^2_{k_2/F}/ d\Omega^1_{k_2/F}$ is 1 and that it is generated by
 the image of $ w=\frac{d b_1}{b_1} \wedge \frac{d b_2}{b_2}$. By assumption, $F$ does not admit any extension of degree prime to $p$ and we have $[k_1:F]=p$ by definition. Hence, the assumptions from Lemma \ref{super-zero} are satisfied (\added[id=Bastian]{note that $g(1)\neq0$ }).
Thus, we obtain an element $z_1 \in k_1 \setminus F$  such that 
\begin{align*}
  z_1^i a \frac{d b_1}{b_1} \wedge \frac{d b_2}{b_2} \in  d \Omega^1_{k_2/F}
\end{align*}
for $1\leq i \leq p-1$. Since  $z_1 \in k_1 \setminus F$, there is an element $a'\in k_2$ such that $a \frac{d b_1}{b_1} \wedge \frac{d b_2}{b_2}= a' \frac{d z_1}{z_1} \wedge \frac{d b_2}{b_2}$ holds in $\Omega^2_{k_2/F}$.
 By Lemma \ref{decomp}, this implies that $a'=e_0 +e_1$ where $e_0 \in F  $ and $e_1\in V$ where $V$ is the $k_1$ vector space generated by
$b_2^i$ for $1\leq i \leq p-1$.
It follows that we can apply Lemma \ref{symbol} to the element $a' \frac{d z_1 }{z_1 } \wedge \frac{d b_2}{b_2}$ which tells us that there is an element
$z_2 \in k_2 \setminus k_1   = F\brac{b_1,b_2}\setminus F\brac{b_1}$ such that $a' \frac{d z_1}{z_1} \wedge \frac{d b_2}{b_2}=\frac{d z_1}{z_1} \wedge \frac{d z_2}{z_2}$ holds in $\Omega^2_{k_2/F}$ and
therefore also in $\Omega^2_{k/F}$. \par

\end{proof}

\begin{lemma}\label{two}
  Let $w=a_1 \frac{d b_1}{b_1} \wedge \frac{d b_2}{b_2}+a_2 \frac{d b_1}{b_1} \wedge \frac{d b_3}{b_3}$ 
be an element in $\nu(2)_k$ with $a_2\neq 0$. Then, there are $z_1 \in F(b_1) \setminus F$ and
$z_2 \in F(b_1,b_2,b_3) \setminus F(b_1,b_2)$ such that 
\begin{align*}
  w= a_1' \frac{d b_1}{b_1} \wedge \frac{d b_2}{b_2}+\frac{d z_1}{z_1} \wedge \frac{d z_2}{z_2}.
\end{align*}
for some $a_1' \in k$.
\end{lemma}
\begin{proof}
 Let us define the fields $k_1=F(b_1), k_2=F(b_1,b_2,b_3)$.
We first claim that $ a_i \in k_2$. As $w = a_2 \frac{d b_1}{b_1} \wedge \frac{d b_3}{b_3}$ 
in $\Omega^2_{k/F(b_2)}$ and $ w \in \nu(2)_k$, it follows that $da_2\wedge \frac{d b_1}{b_1} \wedge \frac{d b_3}{b_3}
=0$ in $\Omega^2_{k/F(b_2)}$. Thus, $a_2 \in k_2$. By symmetry, it follows that $a_1 \in k_2$.
Hence, our assumptions carry over to $\Omega^2_{k_2/F}$. 

Consider the map $s: \{1,2\} \rightarrow \{1,2,3\}\in S_2$
defined via $s(1)=1$ and $s(2)=3$. This just means $\frac{d b_1}{b_1} \wedge \frac{d b_3}{b_3}=w_s$ (compare Subsection \ref{filtration}). Let us furthermore
define $w_{\max}=\frac{d b_1}{b_1} \wedge \frac{d b_3}{b_3}\wedge \frac{db_2}{b_2}$.\par
Note that, by definition of $\gamma-1$, we have $(\gamma-1)(\Omega_{k_2/F,<s})\subset \Omega_{k_2/F,<s}$. 
As $(\gamma-1)(w)$ is a boundary and  $(\gamma-1)(a_1 \frac{d b_1}{b_1} \wedge \frac{d b_2}{b_2}) \in  \Omega^2_{k_2,<s}$, we can  conclude
that 
\begin{align*}
(\gamma-1)\left(a_2 \frac{d b_1}{b_1} \wedge \frac{d b_3}{b_3}\right) \in \Omega^2_{k_2/F,<s}+d\Omega^1_{k_2/F}
\end{align*}
holds.  \par
Consider the map 
\begin{align*}
  g \colon k_1 a_2 w_{\max} \xrightarrow{\textrm{pr}} \Omega^3_{k_2/F} / d\Omega^2_{k_2/F}.
\end{align*}
As $\Omega^3_{k/F} / d \Omega^2_{k/F}$ is a 1-dimensional $F$ vector space \fxnote{Like above, more explanations?}, we can apply Lemma \ref{super-zero} \added[id=Bastian]{(note that $g(1)\neq 0)$)}
to obtain an element $z_1 \in k_1 \setminus F$ such that 
\begin{align}
  z_1^i a_2 w_{\max}\in d \Omega^2_{k_2/F}
  \label{ass:lem}
\end{align}
holds. As $k_1=F(z_1)$, the set $\{z_1,b_2,b_3\}$ forms a $p$-basis of $k_2$ over $F$. Hence, we can rewrite 
\begin{align*}
  a_2 \frac{d b_1}{b_1} \wedge \frac{d b_3}{b_3}=a_2' \frac{d z_1}{z_1} \wedge \frac{d b_3}{b_3} \in \Omega^2_{k_2/F}
\end{align*}
for some $a_2'\in k_2$. But now, taking statement (\ref{ass:lem}) into consideration, Lemma \ref{decomp} guarantees the existence of $e_0 \in F$ and $e_1 \in V$
such that $a_2'=e_0+e_1$ where $V=\bigoplus_{i,j} k_1 b_2^i b_3^j$ where $0\leq i,j \leq p-1$ and $i+j>0$.
We now have 
\begin{align*}
  (a_2'^p-a_2')\frac{d z_1}{z_1} \wedge \frac{d b_3}{b_3} \in \Omega^2_{k_2/F,<s}+d \Omega^1_{k_2/F},
\end{align*}
which by noting that $\frac{db_2}{b_2}\wedge \Omega^2_{k_2/F,<s}=0$ yields 
\begin{align*}
  (a_2'^p-a_2')\frac{d z_1}{z_1} \wedge \frac{d b_3}{b_3} \wedge \frac{db_2}{b_2} \in d \Omega^2_{k_2/F}.
\end{align*}
Therefore, in regard of Lemma \ref{exact}, we can deduce that the zero component, which is given by
\begin{align*}
  (a_2'^p-e_0') \frac{d z_1}{z_1} \wedge \frac{d b_3}{b_3} \wedge \frac{db_2}{b_2},
\end{align*}
is zero. Hence, we obtain $a_2'^p-a_2'=-e_1\in V$. This gives us 
\begin{align*}
  e_1 \frac{dz_1 }{z_1} \wedge \frac{d b_3}{b_3} \in \Omega^2_{k_2,<s}+ d \Omega^1_{k_2/F}, 
\end{align*}
which in turn implies 
\begin{align*}
  d\left( e_1\frac{dz_1}{z_1}\wedge\frac{db_3}{b_3}\right) \in d \Omega^2_{k_2,<s}.
\end{align*}
Recalling $k_1 = F(b_1)=F(z_1)$ now shows us 
\begin{align*}
  \Omega^2_{k_2,<s} = \Omega^1_{k_2,<s(2)}\wedge \frac{dz_1}{z_1},
\end{align*}
which explains the existence of an element $l \in k_2$ such that 
\begin{align*}
  d \left( e_1 \frac{db_3}{b_3}-l \frac{db_2}{b_2} \right) \wedge \frac{dz_1}{z_1}=0 \in \Omega^2_{k_2/F}
\end{align*}
holds. Using Lemma \ref{cycletest}, we get 
\begin{align*}
   d \left( e_1 \frac{db_3}{b_3}-l \frac{db_2}{b_2} \right)=0 \in \Omega^2_{k_2/k_1}.
\end{align*}
We can assume without loss of generality that $l \frac{db_2}{b_2}$ has no zero component over $k_1$, as this would lie
in the kernel of the differential anyway. As, by definition, $e_1' \frac{db_3}{b_3}$ also has no zero component over $k_1$,
Lemma \ref{exact} implies 
\begin{align*}
  e_1 \frac{db_3}{b_3}-l\frac{db_2}{b_2} \in d \Omega^0_{k_2/k_1},
\end{align*}
where we write $\Omega^0_{k_2/k_1}$ instead of $k_2$ to emphasize that this statement is relative to $k_1$.
We can rephrase the last statement as 
\begin{align*}
  (a_2'^p-a_2')\frac{db_3}{b_3}=-e_1\frac{db_3}{b_3} \in d \Omega^0_{k_2/k_1}+\Omega^1_{k_2/k_1,<s(2)}
\end{align*}
which  implies 
\begin{align*}
   (a_2'^p-a_2')\frac{db_3}{b_3} \in d \Omega^0_{k_2/k_1(b_2)}.
\end{align*}
Since the $p$-rank of $k$ over $k_1(b_2)$ is 1, it  follows immediately that $a_2' \frac{db_3}{b_3} \in Z^1_{k_2/k_1(b_2)}$ holds.
Thus, we can apply Theorem \ref{cartier} to obtain $z_2\in k_2 \setminus k_1(b_2)$ such that 
\begin{align*}
  a_2'\frac{db_3}{b_3}=\frac{dz_2}{z_2} \in \Omega^1_{k_2/k_1(b_2)}.
\end{align*}
In view of Lemma \ref{kernel}, this gives us 
\begin{align*}
   a_2'\frac{db_3}{b_3}=\frac{dz_2}{z_2}+l_1\frac{db_1}{b_1} + l_2 \frac{db_2}{b_2} \in \Omega^1_{k_2/F}
\end{align*}
for some $l_i \in k_2$.
So, we eventually get 
\begin{align*}
  a_2 \frac{d b_1}{b_1} \wedge \frac{d b_3}{b_3}=a_2' \frac{d z_1}{z_1} \wedge \frac{d b_3}{b_3} = \frac{d z_1}{z_1} \wedge \frac{d z_2}{z_2} + l_2\frac{d z_1}{z_1} \wedge \frac{d b_2}{b_2}
\end{align*}
in $\Omega^2_{k_2/F}$. \par
By taking into account that there is some $l \in k_2 $ such that $l_2\frac{d z_1}{z_1} \wedge \frac{d b_2}{b_2}=l \frac{d b_1}{b_1} \wedge \frac{d b_2}{b_2}$ holds, we can rewrite our original expression 
\begin{align*}
 & a_1 \frac{d b_1}{b_1} \wedge \frac{d b_2}{b_2}+a_2 \frac{d b_1}{b_1} \wedge \frac{d b_3 }{b_3} \\=&
  a_1' \frac{d b_1}{b_1} \wedge \frac{d b_2}{b_2}+ \frac{d z_1}{z_1} \wedge \frac{d z_2 }{z_2} 
\end{align*}
which holds in $\Omega_{k_2/F}$ and thus in $\Omega_{k/F}$.\par
\end{proof}

\begin{lemma}\label{three}
    Let $w=a_1 \frac{d b_1}{b_1} \wedge \frac{d b_2}{b_2}+a_2 \frac{d b_1}{b_1} \wedge \frac{d b_3}{b_3}+a_3 \frac{d b_2}{b_2} \wedge \frac{d b_3}{b_3}$ 
be an element in $\nu(2)_k$ with $a_3 \neq 0$. Then, there are $z_1 \in F(b_1,b_2) \setminus F(b_1)$ and
$z_2 \in F(b_1,b_2,b_3) \setminus F(b_1,b_2)$ such that 
\begin{align*}
  w= a_1' \frac{d b_1}{b_1} \wedge \frac{d b_2}{b_2}+a_2' \frac{d b_1}{b_1} \wedge \frac{d b_3}{b_3}
 +\frac{d z_1}{z_1} \wedge \frac{d z_2}{z_2}.
\end{align*}
for some $a_1',a_2' \in k$.
\end{lemma}
\begin{proof}
 Let us define the fields $k_0=F(b_1), k_1=F(b_1,b_2)$ and $k_2=F(b_1,b_2,b_3)$. 
Just like in the proof of Lemma \ref{two}, we can deduce that $a_i \in k_2$ and thus our assumptions
carry over to $\Omega^2_{k_2/F}$.\par
Recall that $\nu(2)_k=\ker(\gamma-1)$ is functorial
in field extensions of $k^p$ inside of $k$. Hence, we can deduce
that $(a_3^p-a_3)\frac{d b_2}{b_2} \wedge \frac{d b_3}{b_3} \in d \Omega_{k_2/k_0} $
holds by assumption as $\frac{d b_1}{b_1} \wedge \frac{d b_2}{b_2}= \frac{d b_1}{b_1} \wedge \frac{d b_3}{b_3}=0 \in \Omega^2_{k_2/k_0}$.
 Note also that therefore $w=a_3\frac{d b_2}{b_2} \wedge \frac{d b_3}{b_3}\in \Omega^2_{k_2/k_0}$.\par
 We claim that there is a $z_1\in k_1 \setminus k_0$ such that
 \begin{align}
  z_1^ia_3\frac{d b_2}{b_2} \wedge \frac{d b_3}{b_3} \in d  \Omega_{k_2/k_0}^2
\label{cia} 
\end{align}
for all $1\leq i \leq p-1$.\par
To see this, consider the $k_0$-linear mapping 
\begin{align*}
 g \colon k_1 a_3 \frac{d b_2}{b_2} \wedge \frac{d b_3}{b_3} \subset \Omega^2_{k_2/k_0} \xrightarrow{\textrm{pr}} \Omega^2_{k_2/k_0} / d \Omega^1_{k_2/k_0} 
\end{align*} 
and note that the dimension of the $k_0$ vector space $\Omega^2_{k_2/k_0}/ d\Omega^1_{k_2/k_0}$ is 1 and that it is generated by
 the image of $\frac{d b_2}{b_2} \wedge \frac{d b_3}{b_3}$. We therefore obtain a map $g:k_1\rightarrow k_0$, which satisfies the assumptions from Lemma \ref{super-zero}.
To see this, note that $k_0$ also does not admit any extension of degree prime to $p$ and that $[k_1:k_0]=p$ by definition.\fxnote{More explanation? Let's wait for Parimala
and Suresh's feedback}\par
We thus obtain  $z_1 \in k_1 \setminus k_0$ \added[id=Bastian]{(note that $g(1) \neq 0$)} such that $z_1^i a_3\frac{d b_2}{b_2} \wedge \frac{d b_3}{b_3} \in d \Omega^1_{k_2/k_0}$ holds for $1\leq i \leq p-1$.


We thus have $k_1=k_0(z_1)$, which implies that we can write 
\begin{align*}
  a_3 \frac{d b_2}{b_2} \wedge \frac{d b_3}{b_3}=a_3' \frac{d z_1}{z_1} \wedge \frac{d b_3}{b_3} \in \Omega_{k_2/k_0}
\end{align*}
for some $a_3' \in k$. \par
Note that this implies that  $(\gamma-1)(a_3' \frac{d z_1}{z_1} \wedge \frac{d b_3}{b_3})$ is a boundary in $\Omega^2_{k_2/k_0}$.
From now on, we will look at the $p$-basis $\{z_1, b_3\}$ of $k$ over $k_0$.
We can apply Lemma \ref{decomp} which allows us to write $a_3'=e_0+e_1$
for $e_0\in k_0$ and $e_1 \in V$ where $V$ is the $k_0\brac{z_1}=k_1$ vector space generated by $b_3^{j}$ for $1\leq j\leq p-1$.
This in turn means that the assumptions of Lemma \ref{symbol} are satisfied.
Thus, we get $a_3 \frac{d b_2}{b_2} \wedge \frac{d b_3}{b_3}=a_3' \frac{d z_1}{z_1} \wedge \frac{d b_3}{b_3}=\frac{d z_1}{z_1} \wedge \frac{d z_2}{z_2} \in \Omega^2_{k_2/k_0}$
for some $z_2 \in k_2 \setminus k_1 = k_2\setminus F\brac{b_1,b_2}$. Hence,  there are
$l_i \in k_2$ such that 
\begin{align*}
  a_3 \frac{d b_2}{b_2} \wedge \frac{d b_3}{b_3}=\frac{d z_1}{z_1} \wedge \frac{d z_2}{z_2}
  + l_1 \frac{d b_1}{b_1} \wedge \frac{d b_2}{b_2}
  + l_2 \frac{d b_1}{b_1} \wedge \frac{d b_3}{b_3}.
\end{align*}
 \par
So, by taking Lemma \ref{kernel} into account, we can rewrite our original expression 
\begin{align*}
 & a_1 \frac{d b_1}{b_1} \wedge \frac{d b_2}{b_2}+a_2 \frac{d b_1}{b_1} \wedge \frac{d b_3 }{b_3} +a_3 \frac{d b_2}{b_2} \wedge \frac{db_3 }{b_3}\\=&
  a_1' \frac{d b_1}{b_1} \wedge \frac{d b_2}{b_2}+a_2' \frac{d b_1}{b_1} \wedge \frac{d b_3 }{b_3} + \frac{d z_1}{z_1} \wedge \frac{dz_2 }{z_2}
\end{align*}
 in $\Omega^2_{k_2/k_0}$ for some $a_1',a_2' \in k_2$. \par
\end{proof}

\begin{theorem}
Let $k$ be a field of characteristic $p$. Suppose that $k$ does not admit
any extensions of degree prime to $p$. Let $k$ have $p$-rank $3$. Then, any element $\theta$
in $\nu(2)_k \subset \Omega_k^2$ can be written in the form
\begin{align*}
  \frac{ dx_1}{x_1}\wedge \frac{d x_2}{x_2}+ \frac{d y_1}{y_1}\wedge \frac{d y_2}{y_2}+ \frac{d z_1}{z_1}\wedge \frac{d z_2}{z_2}
\end{align*}
where $\{ x_1,y_2,z_1\}$ is a  $p$-basis of $k$ over $F$.
\label{rep}
\end{theorem}
\begin{proof}
  Let us fix a $p$-basis $\{b_1,b_2,b_3\}$ of $k$ over $k^p$. Then, there are $a_i \in k$ such that
  $\theta=a_1 \frac{d b_1}{b_1} \wedge \frac{d b_2}{b_2} + a_2\frac{d b_1}{b_1} \wedge \frac{d b_3}{b_3}
  + a_3 \frac{d b_2}{b_2} \wedge \frac{d b_3}{b_3}
  $
holds. We will proceed in three steps according to the three previous lemmata. \par

If $a_3 \neq 0$, using Lemma \ref{three}, we can deduce that there are $z_1 \in F(b_1,b_2) \setminus F(b_1)$ and $z_2 \in F(b_1,b_2,b_3) \setminus F(b_1,b_2)$ such that 
\begin{align*}
  & a_1 \frac{d b_1}{b_1} \wedge \frac{d b_2}{b_2}
+a_2 \frac{d b_1}{b_1} \wedge \frac{d b_3}{b_3}
+a_3 \frac{d b_2}{b_2} \wedge \frac{d b_3}{b_3}\\
=& a_1' \frac{d b_1}{b_1} \wedge \frac{d b_2}{b_2}
+a_2' \frac{d b_1}{b_1} \wedge \frac{d b_3}{b_3}+ \frac{d z_1}{z_1}\wedge \frac{d z_2}{z_2}
\end{align*}
for some $a_1',a_2' \in k$. Note that 
\begin{align*}
 \theta'= a_1' \frac{d b_1}{b_1} \wedge \frac{d b_2}{b_2}
+a_2' \frac{d b_1}{b_1} \wedge \frac{d b_3}{b_3} \in \nu(2)_k
\end{align*}
as $\theta'=\theta-\frac{d z_1}{z_1} \wedge \frac{d z_2}{z_2}$ and $\theta, \frac{d z_1}{z_1} \wedge \frac{d z_2}{z_2}\in \nu(2)_k$. \par
If $a_3 =0$, then choose $z_1=b_2$ and $z_2=0$. In this case, set $a_2'=a_2$ and $a_1'=a_1$.

If $a_2'\neq 0$, by applying Lemma \ref{two} on $\theta'$, we can conclude that there are $y_1 \in F(b_1) \setminus F$ and $y_2 \in F(b_1,b_2,b_3) \setminus F(b_1,b_2)$
such that 
\begin{align*}
   & a_1 \frac{d b_1}{b_1} \wedge \frac{d b_2}{b_2}
+a_2 \frac{d b_1}{b_1} \wedge \frac{d b_3}{b_3}
+a_3 \frac{d b_2}{b_2} \wedge \frac{d b_3}{b_3}\\
=& a_1'' \frac{d b_1}{b_1} \wedge \frac{d b_2}{b_2}
+ \frac{d y_1}{y_1} \wedge \frac{d y_2}{y_2}+ \frac{d z_1}{z_1}\wedge \frac{d z_2}{z_2}
\end{align*}
for some $a_1'' \in k$.\par
If $a_2'=0$, then choose $y_1=0,y_2=b_3$ and set $a_1''=a_1'$. 

Again, $\theta''=a_1'' \frac{d b_1}{b_1} \wedge \frac{d b_2}{b_2}$ is in $\nu(2)_k$. Consequently, if $a_1''\neq 0$, Lemma \ref{one} ensures the existence of
$x_1 \in F(b_1) \setminus F$ and $x_2 \in F(b_1,b_2) \setminus F(b_1)$ such that 
\begin{align*}
     & a_1 \frac{d b_1}{b_1} \wedge \frac{d b_2}{b_2}
+a_2 \frac{d b_1}{b_1} \wedge \frac{d b_3}{b_3}
+a_3 \frac{d b_2}{b_2} \wedge \frac{d b_3}{b_3}\\
=&  \frac{d x_1}{x_1} \wedge \frac{d x_2}{x_2}
+ \frac{d y_1}{y_1} \wedge \frac{d y_2}{y_2}+ \frac{d z_1}{z_1}\wedge \frac{d z_2}{z_2}
\end{align*}
holds. \par
If $a_1''=0$, choose $x_1=b_1$ and $x_2=0$.

It is also clear from above that $\{x_1,y_2,z_1\}$ form a $p$-basis. 
\end{proof}

\begin{remark}
Let $k$ be a field of characteristic $p$. Suppose that $k$ does not admit
any extensions of degree prime to $p$. Let $k/F$ have $p$-rank $2$. Then, any nonzero element $\theta$
in $\nu(2)_k \subset \Omega_k^2$ can be written in the form
\begin{align*}
  \frac{ dz_1}{z_1}\wedge \frac{d z_2}{z_2}
\end{align*}
where $\{z_1,z_2\}$ is a  $p$-basis of $k$ over $k^p$. This follows readily from the proof of the surjectivity of $ k_2(k)\rightarrow \nu(2)_{k}$ (compare
also step 3 in the proof of Theorem \ref{rep}). 
  \label{rep_p=2}
\end{remark}

\section{The Brauer $p$-dimension}
Let us recall the definition of Brauer $p$-dimension.
\begin{definition}
  For any field $K$, the \emph{Brauer $p$-dimension} of $K$, denoted by $\textrm{Br}_p \dim (K)$, is the smallest integer $d\geq 0$ such that for any
finite field extension $E$ of $K$ and any central simple algebra $A$ of period a power of $p$ over $E$ we have $\ind (A)|\per (A)^d$.
\end{definition}
 Parimala and Suresh derived a general bound
for the Brauer-$p$-dimension of a complete discretely valued field in terms of the $p$-rank of its residue field. 
\begin{theorem}[\cite{Parimala2013}, Thm 2.7]
  Let $K$ be a complete discretely valued field of characteristic 0  with residue field $\kappa$. Assume that $\kappa$ has characteristic $p$ and $p$-rank $n$ over
$\kappa^p$. If $n=0$, then $\Br_p \dim (K)\leq 1$ and otherwise $\lfloor n/2 \rfloor \leq \Br_p \dim (K)\leq 2n$.
\label{bound}
\end{theorem} 

Note that when determining the Brauer $p$-dimension, we only need to consider 
the case where the period is $p$ (cf. Lemma 1.1 \cite{Parimala2013}).
Let us recall another  well-known reduction.

\begin{remark}
\label{reduc}
  When determining the Brauer $p$-dimension of a field $K$ of characteristic 0, we may assume without loss of generality that $K$ does not admit any extension of degree prime to $p$.
 This follows from the following  argument (cf. Theorem 2.4 in  \cite{Parimala2013}):
 
For otherwise, note that for an extension $K'/K$ of degree prime to $p$, by Lemma \ref{prank}, we have that the $p$-rank of $\kappa'$ equals the $p$-rank of $\kappa$ where $\kappa'$ denotes the residue field of $K'$.
Also, for any central simple algebra $D$ of period $p$ over $K$, we get $\ind(D)=\ind (D \otimes_{K} K')$ as $[K':K]$ is coprime to $p$. 

Note that this implies in particular that we can without loss of generality assume that $K$ contains a primitive $p$-th root of unity. 
\end{remark}

We also recall a filtration of the $p$-torsion part of the Brauer group and its ties to differentials 
and  Milnor $k$-groups.\par
Let $(R,\nu)$ be a complete discrete valuation ring of characteristic 0  with field of fractions $K$ and residue field $\kappa$.
 Let $\pi\in R$ be a parameter. Assume that $K$ has a $p$-th root of unity and let $N=\nu(p)p/(p-1) $. Note that $N$ is a postive integer (cf. Proposition 4.1.2 in \cite{Colliot-Thelene1997}).  
 
Let $\br (K)_0=\prescript{}{p}{\Br(K)} $ and let $U_i=\left \{ u \in R^{*} \; | \; u \equiv 1 \textrm{ mod } \brac{\pi^i} \right \}$. Then, for $i\geq 1$,  let $\br (K)_{i}$ be the subgroup of $\br (K)_0$ generated by cyclic algebras of the form $(u,a)$ for $u \in U_i$ and $a \in K^{*}$. As $R$ is complete, every element in $U_{n}$ is a $p$-th power  whenever $n>N$  so that $\br (K)_n=0$ (cf. Proposition 4.1.2 in \cite{Colliot-Thelene1997}). \par
For arbitrary $z \in \kappa$, let $\tilde{z}$ denote a lift in $R$. Consider the (not functorial!) homomorphism 
\begin{align*}
  \rho_0\colon k_2(\kappa) \oplus \kappa^* / \kappa^{*p} &\rightarrow \br(K)_0 / \br(K)_1, \\
  ((x,y),z)  &\mapsto  (\tilde{x},\tilde{y})+(\pi,\tilde{z}). \\
\end{align*}
Kato proved that $\rho_0$ is in fact an isomorphism (cf. Theorem 2 in \cite{Kato1982} and  Theorem 4.3.1 in \cite{Colliot-Thelene1997}). In the same paper, Kato also proved that the morphism defined via
\begin{align*}
\rho_1\colon \Omega^1_{\kappa} &\rightarrow \br(K)_1 / \br(K)_2\\
  y \frac{dx}{x} &\mapsto (1+\pi \tilde{y},\tilde{x})
\end{align*}
is an (again, not functorial!) isomorphism. The filtration of $\br(K)_0$ is shorter in the case when the residue field is separably closed.
\begin{lemma}[\cite{Colliot-Thelene1997}, Prop 4.2.4]
  Let $\kappa$ be separably closed. Then, $\br(K)_N=0$.
\end{lemma} 
 Recall the following proposition for central simple algebras in $\br(K)_1$.
\begin{proposition}[\cite{Parimala2013}, Proposition 2.2]
  Let $R$, $K$, $\kappa$ and $\pi$ be as above. Suppose that $\kappa=\kappa^p(a_1,\ldots,a_n)$ for some  $a_i \in \kappa$. Let $\alpha \in \br(K)_1$.
Then, for any choice of lifts $\tilde{a}_i$, there are $\lambda, \lambda_i \in R^{*}$ such that 
\begin{align*}
  \alpha=(\lambda_1,\tilde{a}_1)+\ldots+(\lambda_n,\tilde{a}_n)+(\pi,\lambda).
\end{align*}
\label{br1_rep}
\end{proposition}  
 \subsection{Main Theorem}
Let $K$ denote a complete discretely valued field of characteristic $0$ with residue field $\kappa$ of characteristic $p$. We now find some bounds for the Brauer-$p$-dimension of $K$ when the $p$-rank of $\kappa$ is small. Let $R$ denote the valuation ring and $\pi$ a parameter.

\begin{remark}
  \label{redux}
  Note that the condition that $K$ is a complete discretely valued field is stable under finite extensions. Also, the $p$-rank of the residue field
  is preserved under finite extension (c.f. Lemma \ref{prank}). Hence, it is enough to consider only central simple algebras over $K$ when determining
  its Brauer-$p$-dimension. \par
  Further, it is enough to consider central simple algebras of period $p$, compare Lemma 1.1 in \cite{Parimala2013}.
\end{remark}

\begin{proposition} 
  Let $[\kappa:\kappa^p]=p$. Then  $1 \leq \Br_p \dim (K) \leq 2$.
  \label{low_p=1}
\end{proposition} 
\begin{proof}
  The upper bound follows from Theorem \ref{bound}. \par
  For the lower bound, note that $\kappa^*/\kappa^{*p}\neq \{1\}$ as the $p$-rank is 1. Additionally,
  $k_2(\kappa)$ is trivial. Hence, the result follows
from the isomorphism 
\begin{align*}
  \kappa^*/ \kappa^{*p} \longrightarrow \br(K)_0 / \br(K)_1.
\end{align*}
\end{proof}
We now turn to the case where the residue field is separably closed.
\begin{proposition}
   If $[\kappa:\kappa^p]=p$ and $\kappa$ is separably closed, then $ \Br_p \dim (K)=1$.
  \label{ksep_p=1}
\end{proposition}
\begin{proof}
Using Proposition \ref{low_p=1}, it is clear that we only need to show that the Brauer-$p$-dimension
is less than or equal to 1. \par
Let  $\alpha \in \br(K)_0$ be an arbitrary element. Then, by a lemma of Saltman (cf. Lemma 2.8 in \cite{Lieblich2012}), there is a parameter $\pi\in K$ such that
$\alpha \otimes K(\sqrt[\leftroot{-3}\uproot{3}p]{\pi})$ is unramified. Denote the ring of integers of $K(\sqrt[p]{\pi})$ by $R'$. Clearly, the residue
field of $R'$ is still $\kappa$, as the ramification index of the extension is $p$. Since $\kappa=\kappa^{\text{sep}}$, $\Br(\kappa)=0$. Therefore, using $\Br(\kappa) \simeq \Br(R')$ (c.f. Corollary 2.13, \cite{Milne1980}), 
there are no non-trivial unramified algebras. Thus we can conclude that $ K(\sqrt[\leftroot{-3}\uproot{3}p]{\pi})$ splits $\alpha$. Therefore, the Brauer-$p$-dimension is exactly 1. 
\end{proof}

We now turn to the case where the $p$-rank is 2.
\begin{theorem} \label{symbol2} 
Assume that $K$ does not admit extensions prime to $p$.
If $[\kappa:\kappa^p]=p^2$, then every central simple algebra representing an element in $\br(K)_0$ is a tensor product of at most three cyclic algebras of degree $p$.  
\end{theorem}
\begin{proof}
Let $\alpha\in \br(K)_0$ be an arbitrary element. Since $\rho_0 \colon k^2(\kappa) \oplus\kappa^*/\kappa^{*p} \rightarrow \br(K)_0 / \br(K)_1$ is an isomorphism, there are
elements $\gamma \in k^2(\kappa)$ and $c \in \kappa^*$ such that $\alpha = \rho_0\brac{\gamma +[c]} \text{ mod } \br(K)_1$. As the $p$-rank of $\kappa$ is 2,
there are $a,b \in \kappa$ such that $\gamma=\{a,b\}$, compare Remark \ref{rep_p=2}. 
From the isomorphism $k_2(\kappa)\rightarrow \nu(2)_{\kappa} \subset \Omega^2_k$, we can conclude that $\{a,b\}=0$ in $k_2(\kappa)$ or $a$ and $b$ form a $p$-basis.\par
Let us first assume that they form a $p$-basis. Note that we have $\alpha-(\tilde{a},\tilde{b})-(\tilde{c},\pi) \in \br(K)_1$ where $\tilde{x}$ denotes a lift of $x$ for $x\in \kappa$. Using Proposition \ref{br1_rep} for the $p$-basis
$\{a,b\}$, we can see that there are $\lambda,\lambda_i \in R^*$ such that 
\begin{align*}  
  \alpha-(\tilde{a},\tilde{b})-(\tilde{c},\pi)=(\lambda_1,\tilde{a})+(\lambda_2,\tilde{b})+(\pi,\lambda) \in \Br(K).
\end{align*}
Hence, it follows that $\alpha=(\tilde{a},f_1)+(\tilde{b},f_2)+(\pi,f_3)$ for some $f_i \in K^*$.\par
If $\{a,b\}$ do not form a $p$-basis, let $\{e,f\}$ denote some $p$-basis of $\kappa$. Then, using Proposition \ref{br1_rep} again, we see that
\begin{align*}
  \alpha-(\tilde{c},\pi)=(\lambda_1,\tilde{e})+(\lambda_2,\tilde{f})+(\pi,\lambda) \in \Br(K)
\end{align*}
for some $\lambda,\lambda_i \in R^*$. Therefore,  $\alpha=(\tilde{e},f_1)+(\tilde{f},f_2)+(\pi,f_3)$ for some $f_i \in K^*$. \par
\end{proof}

\begin{corollary} 
 If $[\kappa:\kappa^p]=p^2$,  then $\Br_p \dim (K)\leq 3$.
  \label{newbound_p=2}
\end{corollary}
\begin{proof}
  Considering Remark \ref{reduc} and Remark \ref{redux}, the result follows from Theorem \ref{symbol2}.
\end{proof}

Let now $[\kappa:\kappa^p]\geq p^2$ and let $a,b \in R$ be such that their residue classes are $p$-independent.
We now want to prove the existence of $c \in R$ such that $(a,b)\otimes(c,\pi)$ is a division algebra.
This implies in particular that $\Br_p \dim(K) \geq 2$.  

\begin{definition}
  Let $a,b \in R$ be such that $\bar{a},\bar{b}$ are $p$-independent.
  Then, we define
  \begin{align*}
    \sab(R)=
    \left \{ c \in R \; \left| \; \bar{c} \notin \kappa^{p} \cup  \sum_{\substack{ 0 \leq i,j \leq p-1 \\ i+j>0}}\bar{a}^i \bar{b}^j \kappa^p  \right \} \right. .
  \end{align*}
\end{definition}

\begin{lemma} \label{csep}
   Let $L/K$ be a field extension and let $S$ be the valuation ring of $L$.
  Let $l$ denote the residue field of $L$ and assume that $l/\kappa$ is separable.
  We then have
  \begin{align*}
    \sab(R) \subset \sab(S).
  \end{align*}
\end{lemma}
\begin{proof}
   First note that since $l/\kappa$ is separable,  $\bar{a}, \bar{b}$ are still $p$-independent (c.f. Lemma \ref{prank}) in $l$ so that $\sab(S)$ is defined. 
   Let $c \in \sab(R)$. It is clear that $\bar{c} \notin l^p$  since $l/\kappa$ is separable and $\bar{c} \notin \kappa^p$.
  We will now show $\bar{c}$ is not in $ \sum_{\substack{ 0 \leq i,j \leq p-1 \\ i+j>0}} \bar{a}^i \bar{b}^j l^p$. \par
  Let $\{\bar{a}=a_1,\bar{b}=a_2,a_3,\ldots,a_n\}$ denote a $p$-basis of $\kappa$.
  Let $a^v = a_1^{v_1} \cdots a_n^{v_n}$ with $0 \leq v_i \leq p-1$ and let $v=\brac{v_1, v_2, \ldots, v_n}$. Set $\mathcal{V}$ to denote
  the set of all such $v$.
  Then for $\bar{c}$, there is a unique expression
  \begin{align*}
    \bar{c}= \sum_{v \in \mathcal{V}} k_v^p a^v
  \end{align*}
  for some $k_v \in \kappa$. By assumption, we have $k_0 \neq 0$ and $k_v \neq 0$ for at least one $v \neq 0$. 
  As $a_i$ remains a $p$-basis of $l$ (c.f. Lemma \ref{prank} again), the claim follows. 
\end{proof}

If $K$ contains a primitive $p$-th root of unity, and $a,b \in R$ are chosen such that
their residues are $p$-independent, then $(a,b)$ is division. This follows for instance from the
isomorphism $\rho_0$. In this case, $\Lambda=\{ w \in (a,b) \; | \: \text{Nrd}(w) \in R \}$ is the unique
maximal order of $(a,b)$ (c.f. Theorem 12.8 in \cite{Reiner1975}).

\begin{lemma}
  Let $K$ contain a primitive $p$-th root of unity $\omega$ and let $a,b \in R$ be such that
  $\bar{a},\bar{b}$ are $p$-independent. Let $\Lambda$ be the unique maximal order in
  the division algebra $D=(a,b)$. Let $\pi_D \in \Lambda$ denote a parameter of the unique extension
  of the valuation of $K$ to $D$. Then, $\Lambda /\langle \pi_D \rangle = \kappa(\sqrt[p]{\bar{a}},\sqrt[p]{\bar{b}})$
  and $\pi_K$ is a parameter of $D$.
\label{order}
\end{lemma} 
\begin{proof}
  Note first that $\sqrt[p]{a}, \sqrt[p]{b} \in \Lambda$. Hence, $\sqrt[p]{\bar{a}}, \sqrt[p]{\bar{b}} \in \Lambda / \langle \pi_D \rangle$. 
  Also, since $\bar{\omega}=1$ and $\sqrt[p]{a}\sqrt[p]{b}=\omega \sqrt[p]{b}\sqrt[p]{a}$, we see that $\sqrt[p]{\bar{a}}, \sqrt[p]{\bar{b}}$ commute.
  Hence $\kappa(\sqrt[p]{\bar{a}}, \sqrt[p]{\bar{b}}) \subset \Lambda /\langle \pi_D\rangle$. Note that $[\Lambda/\langle \pi_D \rangle: \kappa] \leq p^2$
 as $[\Lambda:R]=[D:K]=p^2$. Since $[\kappa(\sqrt[p]{\bar{a}}, \sqrt[p]{\bar{b}}):\kappa]=p^2$, we obtain $\kappa(\sqrt[p]{\bar{a}}, \sqrt[p]{\bar{b}})=\Lambda/\langle \pi_D \rangle$. 
  This also implies that $\pi_K$ is a parameter of $D$ as the ramification index $e$ of $D$ is 1 (as $f=[\Lambda/\langle\pi_D\rangle : \kappa]=p^2=[\Lambda:R]$). 
\end{proof}

The authors would like to thank David J. Saltman for the idea behind the following proposition: 
\begin{proposition} \label{division}
  Assume that $K$ contains a primitive $p$-th root of unity $\omega$.
  Let $a,b$ be as above and let $c \in \sab(R)$. Then, $(a,b) \otimes_K K(\sqrt[p]{c})$
  is a division algebra.
\end{proposition}
\begin{proof}
  Assume that $(a,b)\otimes_K K(\sqrt[p]{c})$ is not a division algebra, so $K(\sqrt[p]{c})$ is a subfield
   of $(a,b)$.
  Hence, $c$ is a $p$-th power in $(a,b)$. Let $\Lambda$ denote the unique maximal order
  of $(a,b)$ and let $\Lambda_0\subset \Lambda$ denote the subset of elements with reduced trace 0.
  As $\text{Nrd}(\sqrt[p]{c})=c \in R$, it follows that $\sqrt[p]{c} \in \Lambda$.
  Observe also that $\sqrt[p]{c}$ has reduced trace 0, 
  which can be easily deduced from the characteristic polynomial.
  Hence, $c \in \text{Nrd}(\Lambda_0)$. \par
  We will now show that this contradicts $c \in \sab(R)$. Note that by Lemma \ref{order},
  we know that $\Lambda/\langle \pi_D \rangle =\kappa(\sqrt[p]{\bar{a}}, \sqrt[p]{\bar{b}})$.
  Let $x,y \in D$ denote the elements such that $x^p=a$ and $y^p=b$.
 Note  that for any
  $\beta = \sum_{0 \leq i,j \leq p-1}f_{ij}x^iy^j$, the reduced trace is $f_{00}$.
  Let $\rho \colon \Lambda \rightarrow \Lambda /\langle \pi_D \rangle$ denote the natural
  projection. 
Then, we obtain 
  \begin{align*}
    \rho(\Lambda_0) =  \sum_{\substack{ 0 \leq i,j \leq p-1 \\ i+j>0}}\sqrt[p]{\bar{a}}^i \sqrt[p]{\bar{b}}^j \kappa.
  \end{align*}
  Note that  $c=\left(\sum_{0 \leq i,j \leq p-1}f_{ij}x^iy^j\right)^p$ for some $f_{ij} \in R$ with $f_{00}=0$. Thus,
  $\overline{c}=\left(\sum_{0 \leq i,j \leq p-1}\overline{f_{ij}}^p\overline{a}^i\overline{b}^j\right)$ and  we obtain 
  \begin{align*}
    \bar{c} \in 
    \sum_{\substack{ 0 \leq i,j \leq p-1 \\ i+j>0}}\bar{a}^i \bar{b}^j \kappa^p
  \end{align*}
  contradicting our choice of $c$.
\end{proof}

We can now prove that $(a,b) \otimes (c,\pi)$ is a division algebra. 
\begin{theorem}
  Assume that $K$ contains a primitive $p$-th root of unity.
  Let $a,b \in R$ be such that $\bar{a},\bar{b}$ are $p$-independent.
  Let $c \in \sab(R)$. Then, 
  \begin{align*}
    (a,b) \otimes_K (c,\pi)
  \end{align*}
  is a division algebra.
\end{theorem}
\begin{proof}
  Let us assume otherwise. Then, the index of $(a,b)\otimes_K(c,\pi_K)$ divides $p$. Let $F$ be a degree $p$-extension splitting $(a,b) \otimes_K (c,\pi_K)$.
  Then,  there is an extension $L/K$ of degree coprime to $p$ such that $LF=L(\sqrt[p]{d})$ holds
  for some $d \in L$. Since $LF$ splits $(a,b) \otimes_K (c,\pi_K)$, we see that 
  \begin{align*}
    (a,b) +(c,\pi_K)=(d,f) \in \Br(L)
  \end{align*}
  for some $f \in L$ (compare Proposition 2.5.3 in \cite{Gille2006}).
  Let $S$ denote the ring of integers in $L$.
  Write $\pi_K=u \pi_L^e$  where $u$ is a unit in $S$ and  $e$ is an integer coprime to $p$. 
  Note that, without loss of generality, we may assume that $d,f \in S$ with $0 \leq \nu_L(d),\nu_L(f) < p$ where $\nu_L$ denotes the discrete valuation
  of $L$. We obtain
  \begin{align*}
    (a,b)+(c,u)+(c^e,\pi_L)=(d,f) \in \Br(L).
  \end{align*}
  Observe that by Lemma \ref{csep}, $c \in \sab(S)$.
  As $e$ is coprime to $p$, we also have $c^e \in \sab(S)$.

  We split the proof into two cases. 
  Assume first that $d,f$ are units in $S$. Let $l$ denote the residue field of $L$. Then, we see that 
  \begin{align*}
    \rho_0^{-1}((d,f))= (\{\bar{d},\bar{f}\},0) \in k_2(l) \oplus l^*/ l^{*p},
  \end{align*}
  while 
   \begin{align*}
    \rho_0^{-1}((a,b)+(c,\pi))= (\{\bar{a},\bar{b}\}+\{\bar{c},\bar{u}\},\bar{c^e}) \in k_2(l) \oplus l^*/ l^{*p}.
  \end{align*}
  Note that  $\bar{c^e} \notin l^{*p}$ as $c^e \in \sab(S)$, this contradicts $(a,b)+(c,\pi)=(d,f)$. \par
  So, $d$ or $f$ is not a unit. Let us assume without loss of generality that $f$ is not a unit.
  Hence, $L(\sqrt[p]{f})$ has residue field $l$, so  $c \in \sab(S')$ where $S'$ denotes the ring of integers of $L(\sqrt[p]{f})$. Note that by Lemma  $\bar{a}, \bar{b}$ are still $p$-independent over $l$. Since $(a,b)+(c,\pi)=(d,f)$, we conclude that $(a,b)$ is split over $L(\sqrt[p]{f})(\sqrt[p]{c})$ contradicting
  Proposition \ref{division}. 
\end{proof}

\begin{corollary} \label{2low}
  Let $[\kappa:\kappa^p] \geq p^2$. Then, $ \Br_p\dim(K) \geq 2$.
\end{corollary}

Let the $p$-rank of $\kappa$ now be 3.

\begin{theorem}\label{symbol3}
  Assume that $K$ does not admit any  extensions prime to $p$.
 If $[\kappa:\kappa^p]=p^3$, then every central simple algebra representing an element in $\br(K)_0$ is a tensor product of at most 4 cyclic algebras
 of degree $p$.
\end{theorem}
\begin{proof} 
Let $\alpha \in \br(K)_0$. As $\rho_0$ is an isomorphim, we can write 
\begin{align*}
  \rho_0^{-1}([\alpha])=\gamma+\bar{u} \in k_2(\kappa) \oplus \kappa^* / \kappa^{*p}.
\end{align*}
 By Theorem \ref{rep},
we obtain $\gamma=\{z'_1,z_2'\}+\{z_1,z_3'\}+\{z_2,z_3\}$ for some $z_i,z'_i \in \kappa$.
Then, using the definition of $\rho_0$, we obtain
\begin{align*}
  \alpha=(\tilde{z}'_1,\tilde{z}'_2)+(\tilde{z}_1,\tilde{z}'_3)+ (\tilde{z}_2,\tilde{z}_3)+(\pi, \tilde{u})+\alpha',
\end{align*}
where $\alpha' \in \br(K)_1$ and the set  $\{z_2',z_1,z_3 \}$ forms a $p$-basis of $\kappa$.
Using Proposition \ref{br1_rep} with $\{z_2',z_1,z_3 \}$ as our $p$-basis, we can find $\lambda_i, \lambda \in R^{*}$ such that 
\begin{align*}
  \alpha'=(\lambda_1,\tilde{z}'_2)+ (\lambda_2,\tilde{z}_1)+(\lambda_3,\tilde{z}_3)+(\pi, \lambda)
\end{align*}
holds. This implies 
\begin{align*}
   \alpha=(\tilde{z}'_1,\tilde{z}'_2)+(\tilde{z}_1,\tilde{z}'_3)+ (\tilde{z}_2,\tilde{z}_3)+(\pi, \tilde{u})+
   (\lambda_1,\tilde{z}'_2)+ (\lambda_2,\tilde{z}_1)+(\lambda_3,\tilde{z}_3)+(\pi, \lambda)
\end{align*}
which in turn tells us that 
\begin{align*}
  \alpha =(\tilde{z}'_2,f_1)+(\tilde{z}_1,f_2)+(\tilde{z}_3,f_3)+(\pi,f_4)
\end{align*}
 for some $f_i \in K$. 
\end{proof}

\begin{corollary} 
 If $[\kappa:\kappa^p]=p^3$, then $2 \leq \Br_p \dim (K)\leq 4$.
\label{newbound}
\end{corollary}
\begin{proof}
  The lower bound follows from Corollary \ref{2low}.
  For the upper bound, considering Remark \ref{reduc} and Remark \ref{redux}, the result follows from Theorem \ref{symbol3}.
\end{proof}
 
In the case where the $p$-rank is odd, we can slightly improve the lower
bound on the $p$-dimension given in Theorem \ref{bound}. 
\begin{proposition} \label{oddlow}
  If  $[\kappa: \kappa^{p}]\geq 2n+1$, then $\Br_p \dim(K) \geq n+1$. 
\end{proposition} 
\begin{proof}
  Let $\{a_{1},\ldots,a_{2n+1} \}\subset R$ be such that 
  $\{\overline{a_1} , \ldots, \overline{a_{2n+1}} \}$ is $p$-independent.
  We claim that the central simple algebra 
  \begin{align*}
    D=(a_1,a_2)+(a_3,a_4)+ \ldots (a_{2n-1},a_{2n})+(a_{2n+1},\pi)=D_0+(a_{2n+1},\pi)
  \end{align*}
  has index $p^{n+1}$. Parimala and Suresh showed in Lemma 2.6 in \cite{Parimala2013} that
  $D_0$ has index $p^n$. In fact, they showed that if $L/K$ is any extension with residue fields
  $l,k$ such that $[l:\kappa] \leq p^{n-1}$, then $D_0 \otimes_K L$ is not split. \par
  Suppose that there is an extension $L/K$ of degree $p^n$ splitting $D$. Let $S$ denote the valuation
  ring of $L$ and let $l$ denote the residue field. We split the proof into two cases. \par
  Assume first that $\pi$ is a parameter of $S$. Then, by use of the isomorphism
  $\rho_0$, we see that $\overline{a_{2n+1}} \in l^{*p}$ and 
  \begin{align*}
    \frac{d \overline{a_1}}{\overline{a_1}} \wedge \frac{d \overline{a_2}}{\overline{a_2}}+ \ldots + 
    \frac{d \overline{a_{2n-1}}}{\overline{a_{2n-1}}} \wedge \frac{d \overline{a_{2n}}}{\overline{a_{2n}}}=0
    \end{align*}
in $\Omega^2_l$. In particular, $\kappa'=\kappa(\sqrt[p]{\overline{a_{2n+1}}}) \subset l$. Note that $\{ \overline{a_1},\ldots, \overline{a_{2n}} \}$ are still $p$-independent
  over $\kappa'$. However, as $[l:\kappa']<p^n$, an application of Lemma 1.6 in \cite{Parimala2013} with $\lambda_i=1$ shows that then
   \begin{align*}
     \frac{d \overline{a_1}}{\overline{a_1}} \wedge \frac{d \overline{a_2}}{\overline{a_2}}+ \ldots +
   \frac{d \overline{a_{2n-1}}}{\overline{a_{2n-1}}} \wedge \frac{d \overline{a_{2n}}}{\overline{a_{2n}}}\neq 0
   \end{align*}
   in $\Omega^2_l$. This however implies $\alpha=\{\overline{a_1},\overline{a_2}\} + \ldots + \{\overline{a_{2n-1}},\overline{a_{2n}} \} \neq 0 \in k_2(l)$.
   As $\rho_0^{-1}([D])=\alpha + [\overline{a_{2n+1}}]$ (even after base changing to $L$), this contradicts the fact that $D$ is split by $L$. \par
  Let us now assume that $\pi$ is not a parameter in $S$. Then, we have $[l:\kappa]\leq p^{n-1}$.
  Consider now $K'=K(\sqrt[p]{a_{2n+1}})$, $L'=L(\sqrt[p]{a_{2n+1}})$ and let $\kappa', l'$ denote
  their residue fields. Then, $[D\otimes_K K']=[D_0 \otimes_K K']$, so it is enough to show
  that $D_0 \otimes_K K'$ is not split over $L'$. But, this is clear, as  $\{\overline{a_1},\ldots, \overline{a_{2n}}\}$ are still $p$-independent over $\kappa'$ and $[l':\kappa']\leq p^{n-1}$.
\end{proof}

We are now able to state the main theorem.
\begin{theorem}
\label{BN-mainthm1}
Let $K$ be a complete discretely valued field with residue field $\kappa$. Suppose that $\chara(\kappa) = p > 0$ and the $p$-rank of $\kappa$ is $n$ where $n=0, 1, 2$ or $3$. Then $\Br_p\dim(K)\leq n+1$. \par
For $n < 3$, we have $n \leq \Br_p \dim(K)$ and for $n=3$ we have $2\leq \Br_p\dim(K)$.
\end{theorem} 
\begin{proof}
  Let us first discuss the upper bound.
  The cases $n=0,1$ follow from Theorem 2.7 (c.f Theorem \ref{bound}) in \cite{Parimala2013}. The cases $n=2,3$ follow from Corollaries \ref{newbound_p=2} and \ref{newbound} respectively. \par
  For $n=0$, the lower bound is trivial. For $n=1$, it follows from Proposition \ref{low_p=1},
  the case $n=2$ was handled in Corollary \ref{2low} and the case $n=3$ follows from Proposition \ref{oddlow}.
\end{proof}

\section{Examples}
\label{examples}
In this section, for each $n\geq 1$, we give an example of a field with residue field of $p$-rank $n$ whose Brauer $p$-dimension is at least $n+1$. We also construct examples realizing all possible Brauer $p$-dimensions as stated in Theorem \ref{BN-mainthm1} for $n\leq 2$.

\subsection{A family of examples}
\label{fexamples}

Fix a prime $p$ and let $k_0=\F_p$. For $n\geq 1$,  set $k_n=\F_p(x_1,\ldots,x_n)$, $E_n = \frac{k_n[t]}{\brac{t^p-t-x_n}}$ and $K_{n} = k_{n-1}((x_{n}))$. Note that the $p$-rank of $k_n$ is $n$. 

The following proposition constructs division algebras over $k_{n}$ which still remain division over $E_{n}$. 

\begin{proposition}
Let $n\geq 1$ and let $D_{n}=[1,x_1)\otimes \ldots \otimes [x_{n-2},x_{n-1})\otimes [x_{n-1},x_{n})$ in $\Br\brac{k_{n}}$. Then, 

\begin{enumerate}
\item
$D_{n}$ has index $p^{n}$ over $K_{n}$  and hence over $k_{n}$
\item
$D_{n}$ has index $p^{n}$ over $E_{n}$
\end{enumerate}

\label{prankn}
\end{proposition} 
\begin{proof}
We will prove the result by induction on $n$.

Let $n=1$. We would like to first show that $D_1 = [1,x_1)$ has index $p$ over $K_1$. Note that $K_1:=\mathbb{F}_p((x_1))$ is a complete discretely valued field with parameter $x_1$. Set $\widetilde{E_0}=K_1[t]/(t^p-t-1)$. Since $t^p-t-1$ is irreducible over $\Fp$, $\widetilde{E_0}$ is an unramified non-split extension of $K_1$. As $x_1$ is a  parameter of $K_1$ it can not be a norm from $\widetilde{E_0}$ by Lemma \ref{norm}. Therefore, by Proposition \ref{normsplit}, $D_1=[1,x_1)$ is nontrivial in $\Br\left(K_1\right)$ and therefore has index $p$. Since $k_1\subseteq K_1$, $D_1$ has index $p$ over $k_1$. 

We claim that $D_1\otimes_{k_1} E_1$ still has index $p$. This is because of the following: Hensel's Lemma implies that  $t^p-t-x_1$ splits in $K_1=\Fp((x_1))$. Thus, we have 
\[E_1\otimes_{k_1} K_1 := K_1[t]/(t^p-t-x_1) = \prod K_1.\] Thus as $D_1$ has index $p$ over $K_1$, it has index $p$ over $E_1\otimes_{k_1}K_1$ and hence over $E_1$. This finishes the base case of our induction. 

Let now $n>1$ and assume that the statement is true for all $r\leq n-1$. Observe that $D_{n}= D_{n-1} \otimes [x_{n-1},x_{n})$. Since $K_{n}$ is a complete discretely valued field with residue field $k_{n-1}$, the algebra $D_{n-1}$ is unramified in $\Br(K_{n})$ (c.f. Lemma \ref{unramified}). 

By the induction hypotheses, we know that $D_{n-1}$ has index $p^{n-1}$ (and hence is division) over $K_{n-1} := k_{n-2}((x_{n-1}))$, over $k_{n-1}$ and over $E_{n-1} :=k_{n-1}[t]/(t^p-t-x_{n-1})$. 

Set $\widetilde{E_{n-1}}:= K_{n}[t]/(t^p-t-x_{n-1})$. Then $\widetilde{E_{n-1}}/K_{n}$ is an unramified non-split cyclic extension and $D_{n-1}$ remains division over $\widetilde{E_{n-1}}$. 
Theorem 5.15 in \cite{Jacob1990}  then immediately implies that $D_n = D_{n-1}\otimes [x_{n-1},x_n)$ is division over $K_n$ and therefore has index $p^n$. Since $k_n \subseteq K_n$, $D_n$ has index $p^n$ over $k_n$ also.

We claim that $D_n\otimes_{k_n} E_n$ still has index $p^n$. The proof is similar to that of the base case. Note that $E_n\otimes_{k_n} K_n := K_n[t]/(t^p-t-x_{n}) = \prod K_n$ and that $D_n$ has index $p^n$ over $K_n$. Therefore it has index $p^n$ over $E_n\otimes_{k_n}K_n$ and hence over $E_n$. \end{proof}

Now we are ready to construct our family of examples. 

\begin{theorem}
  Let $K$ be a complete discretely valued field of characteristic $0$ with parameter $\pi$ and residue
  field $k_n=\mathbb{F}_p(x_1,\ldots, x_n)$ where $n\geq 0$. Let $\tx_i$ denote a lift of $x_i$. Then, the algebra 
  \begin{align*}
    A=[1,\tx_1)\otimes [\tx_1, \tx_2) \otimes  \ldots \otimes [\tx_{n-1},\tx_n)\otimes [\tx_n,\pi)
  \end{align*}
  has index $p^{n+1}$ over $K$.
  \label{n+1char0}
\end{theorem} 
\begin{proof}
If $n=0$, then $A=[1,\pi)$ and the theorem follows from Theorem 5.15 in \cite{Jacob1990} with $I=K$.
For $n \geq 1$, note that $A'=[1,\tx_1)\otimes \ldots \otimes  [\tx_{n-1},\tx_{n})$ is unramified over $K$ and that the residue algebra has index $p^n$ over $E_n$ by Proposition \ref{prankn}. Hence $A'$ has index $p^n$ over the unramified extension $K[t]/(t^p-t-\tx_n)$. Another application of Theorem 5.15 in  \cite{Jacob1990} shows that $A=A' \otimes [\tx_n,\pi)$ has index $p^{n+1}$.

 


\end{proof}

\begin{remark} Theorem \ref{n+1char0} shows that for $n\leq 3$, the upper bounds for Brauer-$p$-dimensions in Theorem \ref{BN-mainthm1} are optimal. This also shows that the optimal upper bound for Brauer-$p$-dimension for a general $n$ is at least $n+1$. 
\end{remark}

\subsection{Examples for $p$-rank $\leq 2$}
We now set forth examples of fields realizing every Brauer $p$-dimension possible for $n\leq 2$. 

\begin{theorem}\label{exampleleq2}
Let $p$ be a prime and let $n, i$ be integers such that $0\leq n\leq 2$ and $n\leq i\leq n+1$. 
Then there exists a characteristic $0$ complete discretely valued field $K$ with residue field $k$ of characteristic $p$ and $p$-rank $n$ whose Brauer-$p$-dimension is $i$.
\end{theorem}

\begin{proof} In this proof, $k$ will always denote a characteristic $p$ field of $p$-rank $n$. Note that there exists a completely discretely valued field $K$ of characteristic $0$ with residue field $k$ (c.f. Thm 2, \cite{MacLane1939}). Let $R$ denote the ring of integers of $K$ and let $\pi\in R$ be a parameter. Finally for any $\theta\in k$, $\tilde{\theta}\in R$ denotes some lift of $\theta$. We break the proof into cases depending on the values of $n$ and $i$. 

\blue{Case I: $i=n+1$} \\
Let $k=\F_p(x_1, \ldots, x_n)$. Theorem \ref{n+1char0} tells us that $\Br_p\dim(K)\geq n+1$. This, in conjuction with Theorem \ref{BN-mainthm1} which shows $\Br_p\dim(K)\leq n+1$, finishes the proof in this case.

\blue{Case II : $i=n=0$} \\
Let $k=\overline{\mathbb{F}_p}$, an algebraic closure of the finite field $\mathbb{F}_p$. Since $k/\mathbb{F}_p$ is separable, the $p$-rank of $k$ is $0$ (c.f. Lemma \ref{prank}). Denote $N=\frac{\nu(p) p}{p-1}$, where $\nu$ is the valuation on $K$. As $k$ is separably closed, by (Thm 4.2.4 in \cite{Colliot-Thelene1997}) we have that $\br(K)_N=0$. As the $p$-rank of $k=0$, it follows that  $k_2(k)=0$ and $\Omega^q_{k}=0$ for $q\geq 1$. Hence, Theorem 4.3.1 in \cite{Colliot-Thelene1997} implies that $\br(K)_0=\br(K)_1 = \ldots =\br(K)_N=0$, which shows that $\Br_p\dim(K)=0$.

\blue{Case III: $i=n=1$}\\
Let $k$ denote a separable closure of $\mathbb{F}_p(x)$. As $\mathbb{F}_p(x)$ has $p$-rank 1, the same holds true for $k$ (c.f Lemma \ref{prank}). It now follows from Proposition \ref{ksep_p=1} that the Brauer $p$-dimension of $K$ is 1.

\blue{Case IV: $i=n=2$}\\
Let $F$ be a complete discretely valued field with residue field $\F_p(x,y)$ and parameter $p$ (consider for instance the fraction field of the completion of $\Z[x,y]_{(p)}$). Let $k$ denote the separable closure of $\F_p(x,y)$. By (Theorem 1, \cite{MacLane1939}), there is an unramified extension of complete discretely valued fields $K' /F$ such
that the residue field of $K'$ is $k$. Let $K=K'(\zeta)$ where $\zeta$ denotes a primitive $p$-th root of unity. Then, the extension $K/K'$ is totally ramified.
Hence, the residue field of $K$ is $k$ and the valuation of $p$ is $p-1$.

 Recall that our bounds give  $\replaced[id=Bastian]{2}{1}\leq \Br_p\dim(K)\leq 3$ (c.f. Corollary \ref{2low}). Hence it suffices to show $\Br_p \dim(K) \leq 2$.

Let $D$ be a central simple algebra. Let $a,b,c \in K$ such that 
\begin{align*}
    \rho_0\colon&k_2(k) \oplus k^* / k^{*p} \rightarrow \br(K)_0 / \br(K)_1 \\
 & (\{a,b\},c)\mapsto  [D] \\
\end{align*}
where $[D]$ denotes the equivalence class of $D$ in $\br(K)_0 / \br(K)_1$. Thus $D-(\tilde{a},\tilde{b})-(\tilde{c},\pi)$ is an element in $\br(K)_1$.  

Recall that the homomorphism
\begin{align*}
    \rho_i \colon&\Omega^i_{k} \rightarrow \br(K)_i / \br(K)_{i+1} \\
 & y\frac{d x}{x} \mapsto (1+\tilde{y}\pi^i,\tilde{x})
 \\
\end{align*}
is an isomorphism for $1\leq i \leq p-1$(cf. Thm 4.3.1 in \cite{Colliot-Thelene1997}). Since we have
$v(p)_K=p-1$, we also have $\br(K)_p=0$ (cf. Proposition 4.2.4 in \cite{Colliot-Thelene1997}). Thus if $\{x,y\}$ form a  $p$-basis of $k$, then any element in $\br(K)_{p-1}$ can be written in the form  $(-,\tilde{x})+(-,\tilde{y})$. Now repeated usage of $\rho_i$ for $1\leq i\leq p-2$ shows that any element in $\br(K)_1$ can also be written in the form $(-,\tilde{x})+(-,\tilde{y})$. 

Let us now prove that there is an extension of degree $p^2$ splitting $D$. For this, we consider
two different cases.

Let us first assume that $a$ and $b$ are $p$-independent, i.e. form a $p$-basis.
If $a$ and $c$ are $p$-independent, then as we have seen before, any element in $\br(K)_1$ can also be written in the form $(-,\tilde{a})+(-,\tilde{c})$. Then $K(\sqrt{\tilde{a}},\sqrt{\tilde{c}})$ splits $D$. If $b$ and $c$ are $p$-independent, then a similar argument shows that $K(\sqrt{\tilde{b}},\sqrt{\tilde{c}})$ splits $D$. 

According to Lemma \ref{ppower}, if neither $\{a,c\}$ nor $\{b,c\}$ is a $p$-basis, then  $c$ is a $p$-th power. Therefore $(\tilde{c},\pi)\in \br(K)_1$. Therefore, $D-(\tilde{a},\tilde{b}) \in \br(K)_1$. Hence by using again that $\{a,b\}$ is a $p$-basis of $k$, it follows that $D-(\tilde{a},\tilde{b})=(u,\tilde{a})+(v,\tilde{b})$ for some $u,v\in K$. Consequently, we can see that $K(\sqrt{\tilde{a}},\sqrt{\tilde{b}})$ splits $D$.

Now,  let us assume that $a$ and $b$ are $p$-dependent. Thus $D-(\tilde{c},\pi)\in \br(K)_1$.
If $c$ is a $p$-th power, then $D\in \br(K)_1$. Thus $D=(u,\tilde{e})+(v,\tilde{f})$ for some $p$-basis $\{e,f\}$ and some $u,v\in K$ so that $K(\sqrt{\tilde{e}},\sqrt{\tilde{f}})$ is a splitting field of $D$. If $c$ is not a $p$-th power, then there is some element $e\in k$ such that $\{c,e\}$ form a $p$-basis. Hence, there are $u,v \in K$ such that $D-(\tilde{c},\pi)=(u,\tilde{c})+(v,\tilde{e})$. Therefore, $K(\sqrt{\tilde{c}},\sqrt{\tilde{e}})$ splits $D$ and we have proven that the
Brauer $p$-dimension of $K$ is at most 2. \end{proof}

\begin{remark} Theorem \ref{exampleleq2} shows that for $n\leq 2$, the lower bounds for Brauer-$p$-dimensions in Theorem \ref{BN-mainthm1} are optimal. \end{remark}

With the support of the low-dimensional examples in this section and Theorem \ref{BN-mainthm1}, we end with the following conjecture:

\begin{conjecture}
\label{conjBH}
Let $K$ be a complete discretely valued field with residue field $\kappa$. Suppose that $\chara(\kappa) = p > 0$ and that the $p$-rank of $\kappa$ is $n$. Then $n \leq \Br_p\dim(K)\leq n+1$. 
\end{conjecture}


\nocite{Cartier1958}
\bibliographystyle{alpha}
\bibliography{library}
\end{document}